\documentclass[oneside]{article}
\usepackage{fontspec}
\let\RussianRomanFamily\rmdefault

\usepackage[russian,english]{babel}
\usepackage{iftex}
\ifPDFTeX\usepackage[utf8]{inputenc}\fi
\usepackage[a4paper,top=3cm,textwidth=137mm,lines=58,headheight=17pt,headsep=9pt,hcentering]{geometry}
\usepackage{amsmath,amssymb,amsthm,graphicx,tikz,caption,url,placeins}
\makeatletter
\let\@afterindentfalse\@afterindenttrue
\@afterindenttrue
\makeatother
\usetikzlibrary{arrows.meta}
\newcommand{\theoremname}{Theorem}
\newtheorem{theorem}{\theoremname}
\newcommand{\lemmaname}{Lemma}
\newtheorem{lemma}{\lemmaname}
\newcommand{\claimname}{Proposition}
\newtheorem{claim}{\claimname}
\theoremstyle{definition}\newcommand{\defnname}{Definition}
\newtheorem{defn}{\defnname}
\newcommand{\F}{\mathbb F_2}
\newcommand{\rank}{\operatorname{rank}}
\newcommand{\rad}{\operatorname{rad}}
\newcommand{\mat}[1]{\begin{pmatrix}#1\end{pmatrix}}
\graphicspath{{images/}}
\title{\textbf{Embeddability of hieroglyphs in the torus and the intersection form}}
\author{Tim Berezin\thanks{Most of this work was carried out at the Faculty of Mathematics and Mechanics of St. Petersburg University between 2018 and 2019. Email: timurovadia@gmail.com}}
\date{October 21, 2026}

\begin{document}
\begingroup
\selectlanguage{english}
\normalfont
\maketitle
\begin{abstract}
We characterize the hieroglyphs that embed in the torus, where a hieroglyph is a one-vertex graph with a rotation system. We prove equivalent criteria in terms of four forbidden hieroglyphs, the structure of the loop graph, reductions to three terminal forms, and the rank of the interlacement matrix. The proof uses the intersection form on the first homology of a disc with bands over $\F$ and the equality between its rank and twice the genus. We describe the effect of reductions and the role of boundary components. Counting the cycles of the boundary permutation yields a linear-time recognition algorithm for an indexed representation of the input.
\end{abstract}
\section*{Introduction}
Let $G$ be a graph with one vertex and finitely many loops, say $n$. A \emph{hieroglyph} is the graph $G$ equipped with a prescribed cyclic order of the loop ends at its unique vertex. Equivalently, it is a one-vertex graph with a rotation system.
\begin{defn}
A hieroglyph is \emph{embeddable in the torus}, or \emph{realizable on the torus}, if its loops can be drawn on the torus without self-intersections or mutual intersections away from their common vertex, with the prescribed cyclic order of the incident ends around that vertex.
\end{defn}
Thus loops may share their vertex but no other points. The embedding need not be cellular. For extensions to graphs with several vertices and to other surfaces, see \cite{en:MT,en:EM}. For a general connected graph, contraction of a chosen spanning tree, with the rotation system carried along, gives a one-vertex graph; homotopy equivalence of the underlying graphs alone is insufficient to preserve the embedding problem.

An \emph{alphabet} is a finite set; its elements are called \emph{letters}. Label the loops by letters and read their ends around the vertex in the prescribed direction. This represents a hieroglyph by a word of length $2n$ in which each letter occurs twice. Words are considered up to relabelling and cyclic shift. The direction is fixed: reversal is not included in this equivalence. The empty hieroglyph is denoted by $()$. A hieroglyph \emph{contains} another if the latter can be obtained by deleting some letters, each together with both its occurrences.

The \emph{loop graph} $H^*$ of a hieroglyph $H$ is the simple graph whose vertices are the loops of $H$. Two vertices are adjacent if their letters alternate, forming $(abab)$. We then say that the loops \emph{interlace}. This is a property of the cyclic order, not an additional intersection point in an embedding.

A \emph{reduction} is a composition of the following operations:
\begin{align*}
(\alpha)\quad &(aaU)\longmapsto(U),\\
(\beta)\quad &(abUbaV)\longmapsto(aUaV),\qquad a\ne b.
\end{align*}
Here $U,V$ are possibly empty subwords, containing no $a$ in $(\alpha)$ and no $a,b$ in $(\beta)$. Each other letter occurs twice in total; in $(\beta)$ its occurrences may lie in the same subword or in different subwords. The two adjacent pairs in $(\beta)$ use four distinct occurrences. A loop removed by $(\alpha)$ is called an \emph{isolated loop}. An isolated vertex of the loop graph need not correspond to such a loop: in $(abbacc)$ the ends of $a$ are not adjacent. This distinction matters in the reduction argument.

A graph is \emph{complete $k$-partite} if its vertex set is partitioned into $k$ independent parts and every two vertices in distinct parts are adjacent. Empty parts are permitted. Throughout, unions of graphs are disjoint unions; the set of isolated vertices may be empty.
\begin{theorem}\label{en:thm:main}
For a hieroglyph $H$, the following conditions are equivalent:
\begin{enumerate}
\item[(A)] $H$ is realizable on the torus;
\item[(B)] $H$ contains none of $(ababcdcd)$, $(abcdabcd)$, $(abacdcbd)$, $(abcadbdc)$;
\item[(C)] $H^*$ is the union of isolated vertices and a complete bipartite or tripartite graph;
\item[(D)] $H$ reduces to $()$, $(abab)$ or $(abcabc)$;
\item[(E)] the interlacement matrix $M(H)$ over $\F$ has rank at most $2$.
\end{enumerate}
\end{theorem}
The equivalence of (A) and (E) follows from $\rank M(H)=2g(H)$ over $\F$, where $g(H)$ is the genus of the surface obtained by capping the boundary components of a disc with bands. We prove this equality in Section~\ref{en:sec:form}. Related embedding criteria are discussed in \cite{en:S}; for the relation between band surfaces and words in free groups, see \cite{en:C}.

\section{A disc with bands and the intersection form}\label{en:sec:form}
Mark $2n$ pairwise disjoint intervals on the oriented boundary of a disc, in the order of the letters of the word. For each letter, attach a rectangular band along its two intervals. Choose the attaching maps so that the orientation of the disc extends over every band; twisted bands are excluded. Denote the resulting connected compact orientable surface by $S(H)$, a \emph{disc with $n$ bands}. Write $F(H)$ for the number of its boundary components. Set $S(())=D^2$ and $F(())=1$. When convenient, we identify a hieroglyph with its disc with bands. Figure~\ref{en:fig:bands} shows $(abcabc)$ and its two boundary walks; the occurrences on the right are numbered $0,\ldots,5$ in the order $a,b,c,a,b,c$.

\begin{figure}[htbp]\centering
\begin{minipage}[c]{.22\linewidth}\centering
\includegraphics[width=.95\linewidth]{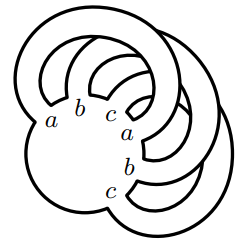}
\end{minipage}\hspace{.06\linewidth}\begin{minipage}[c]{.46\linewidth}\centering
\begin{tikzpicture}[scale=.8,>=Stealth]
\foreach \x/\col/\title/\top/\left/\right in {0/red/{c_1:(0\ 4\ 2)}/{0:a}/{2:c}/{4:b},3.4/blue/{c_2:(1\ 5\ 3)}/{1:b}/{3:a}/{5:c}}{
\begin{scope}[shift={(\x,0)},color=\col]
\draw[->] (90:.75) arc (90:-30:.75);\draw[->] (-30:.75) arc (-30:-150:.75);\draw[->] (-150:.75) arc (-150:-270:.75);
\node at (0,1.55) {$\title$};\node at (0,1.02) {$\top$};\node at (-1.05,-.68) {$\left$};\node at (1.05,-.68) {$\right$};
\end{scope}}
\end{tikzpicture}\end{minipage}
\caption{The hieroglyph $(abcabc)$: a disc with three bands (left) and its two boundary walks (right). The label $i:x$ denotes the side of band $x$ associated with occurrence $i$; arrows indicate the boundary traversal.}\label{en:fig:bands}
\end{figure}

Contracting the disc to the vertex and the bands to their centre lines gives a deformation retraction of $S(H)$ onto $H$. In particular,
\[
\chi(S(H))=1-n,\qquad H_1(S(H);\F)\cong\F^n.
\]
Cap every boundary component by a disc to obtain a closed orientable surface $\widehat S(H)$, whose genus is denoted by $g(H)$.
\begin{lemma}\label{en:lem:genus}
For every hieroglyph,
\begin{equation}\label{en:eq:euler}2g(H)=n+1-F(H).
\end{equation}
It is realizable on a closed orientable surface of genus $h$ if and only if $g(H)\le h$.
\end{lemma}
\begin{proof}
Capping the $F(H)$ boundary components increases the Euler characteristic by $F(H)$, so $2-2g(H)=1-n+F(H)$. If $H$ is embedded in a surface of genus $h$, a closed regular neighbourhood is homeomorphic to $S(H)$. Suppose the closure of its complement has $k$ components, of genera $h_1,\ldots,h_k$, with $F(H)$ boundary components in total. Additivity of $\chi$ under gluing along circles gives
\[
h=g(H)+\sum_{i=1}^k h_i+F(H)-k\ge g(H),
\]
since each complementary component has nonempty boundary. Conversely, embed $H$ in $\widehat S(H)$ and, if necessary, add $h-g(H)$ handles inside a complementary disc. This preserves the rotation system.
\end{proof}
Thus the number of complementary components of an arbitrary embedding must not be identified with $F(H)$: one component can have several boundary circles. Formula~\eqref{en:eq:euler} concerns the canonical capped surface, in which the embedding is cellular.
\begin{defn}
The \emph{interlacement matrix} $M(H)=(m_{ij})_{i,j=1}^n$ over $\F$ has zero diagonal and $m_{ij}=1$ precisely when loops $i,j$ interlace. It is the adjacency matrix of $H^*$ over $\F$.
\end{defn}
For the intersection form $\lambda$ on $H_1(S(H);\F)$, its \emph{radical} is
\[
\rad\lambda=\{x:\lambda(x,y)=0\text{ for all }y\in H_1(S(H);\F)\}.
\]
In the loop basis this radical is identified with $\ker M(H)$; we also write $\rad M(H)$.
\begin{lemma}\label{en:lem:rank}
The matrix $M(H)$ represents the intersection form on $H_1(S(H);\F)$ in the loop basis. Its radical is generated by the boundary classes, and the only relation among these classes is that their sum is zero. In particular,
\begin{equation}\label{en:eq:rank}
\dim\rad M(H)=F(H)-1,\qquad \rank M(H)=2g(H).
\end{equation}
\end{lemma}
\begin{proof}
Represent each basis cycle by the centre line of its band and an arc in the disc joining the ends. Choose the arcs in general position. Their intersection number modulo $2$ is $1$ exactly when their ends alternate on the boundary; cycles in distinct bands do not meet outside the disc. Self-intersection is zero because the surface is orientable. This proves the matrix assertion.

For completeness, consider the standard basis for an orientable surface of genus $g=g(H)$ with $F=F(H)$ boundary components. Choose handle cycles $u_1,v_1,\ldots,u_g,v_g$ and boundary cycles $c_1,\ldots,c_{F-1}$. Cutting along handle cycles and arcs between boundary components gives the standard polygon and shows that these $2g+F-1=n$ classes form a basis. The only nonzero pairwise intersections are $u_i\cdot v_i=1$. The matrix of the form is therefore
\[
\bigoplus_{i=1}^g\mat{0&1\\1&0}\ \oplus\ 0_{F-1}.
\]
The last boundary component is homologous to the sum of the others. This describes the radical and proves~\eqref{en:eq:rank}.
\end{proof}
There is also a natural homological formulation. Write $S=S(H)$ and $\widehat S=\widehat S(H)$. For the inclusion $j:S\hookrightarrow\widehat S$, the exact sequence of the pair and excision of the capping discs give, throughout with coefficients in $\F$,
\[
0\longrightarrow H_2(\widehat S)\longrightarrow H_2(\widehat S,S)\longrightarrow H_1(S)\xrightarrow{j_*}H_1(\widehat S)\longrightarrow0.
\]
The first two nonzero groups are $\F$ and $\F^F$. The first nonzero map sends $1$ to $(1,\ldots,1)$, and the next sends the relative homology classes of the capping discs to their boundary classes. Consequently,
\[
H_1(S(H);\F)/\rad M(H)\cong H_1(\widehat S(H);\F).
\]
This is an isomorphism of spaces with intersection forms; the form on the right is nondegenerate. Thus the radical records the boundary, while the nonzero part of the form records the handles. For standard facts about exact sequences of pairs and intersection forms on orientable surfaces, see \cite{en:H}.

Return to $(abcabc)$ in Figure~\ref{en:fig:bands}. Each of its two boundary walks uses one side of each band. After retraction onto the graph, both walks therefore give the same class:
\[
[c_1]=[c_2]=e_a+e_b+e_c\quad\text{in }H_1(S(H);\F).
\]
Here $e_a,e_b,e_c$ are the loop basis classes. Indeed,
\[
M(abcabc)=\mat{0&1&1\\1&0&1\\1&1&0},\qquad
M(abcabc)\mat{1\\1\\1}=0.
\]
The radical is one-dimensional and generated by this nonzero class; the unique independent relation among the boundary classes is $[c_1]+[c_2]=0$. Lemma~\ref{en:lem:boundary} justifies the two cycles displayed in the figure. Figure~\ref{en:fig:intersection} illustrates an individual entry of the form: alternating ends give one intersection modulo $2$ in the disc.

\begin{figure}[htbp]\centering
\begin{tikzpicture}[scale=.75]
\draw (0,0) circle (1);\draw[thick] (-.707,.707)--(.707,-.707);\draw[thick] (.707,.707)--(-.707,-.707);
\node at (-.9,.9) {$a$};\node at (.9,.9) {$b$};\node at (.9,-.9) {$a$};\node at (-.9,-.9) {$b$};
\node[right] at (2,0) {$\displaystyle M(abab)=\mat{0&1\\1&0},\quad\rank M=2,\ g=1,\ F=1.$};
\end{tikzpicture}
\caption{Alternating ends give an entry $1$ in the intersection form. The crossing of arcs inside the disc computes the form on cycles; it is not a crossing of loops in an embedding.}\label{en:fig:intersection}
\end{figure}

\section{Proof of Theorem~\ref{en:thm:main}}
\subsection*{Forbidden hieroglyphs and small matrices}
Figure~\ref{en:fig:forbidden} shows the four forbidden hieroglyphs in (B). Their loop graphs in Figure~\ref{en:fig:graphs} determine their interlacement matrices: an edge between distinct letters gives an entry $1$.

\begin{figure}[htbp]\centering\includegraphics[width=\linewidth]{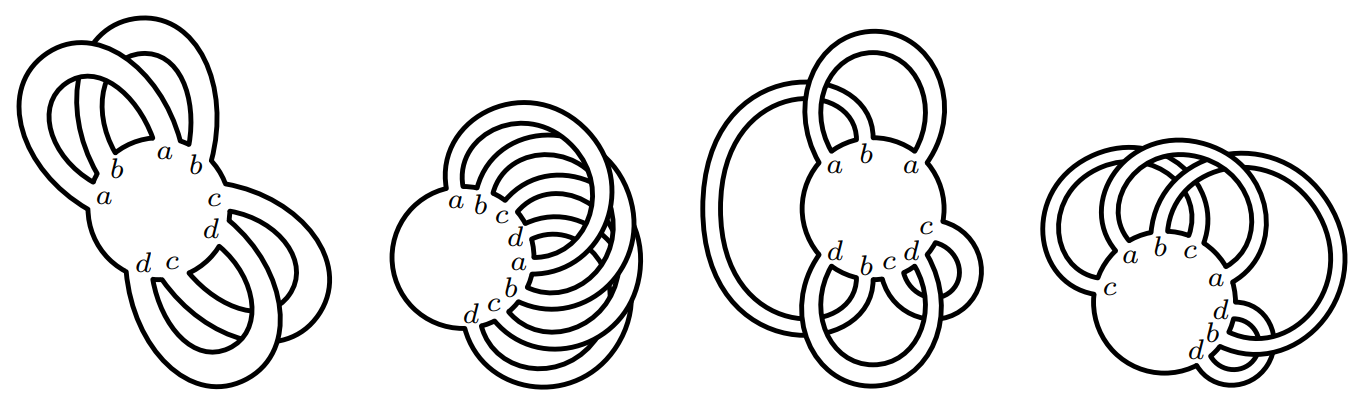}
\caption{The hieroglyphs $(ababcdcd)$, $(abcdabcd)$, $(abacdcbd)$, $(abcadbdc)$.}\label{en:fig:forbidden}\end{figure}

In the listed order and in the basis $a,b,c,d$, these matrices are
\begin{align*}
M_1&=\mat{0&1&0&0\\1&0&0&0\\0&0&0&1\\0&0&1&0},&
M_2&=\mat{0&1&1&1\\1&0&1&1\\1&1&0&1\\1&1&1&0},\\
M_3&=\mat{0&1&0&0\\1&0&0&1\\0&0&0&1\\0&1&1&0},&
M_4&=\mat{0&1&1&0\\1&0&1&1\\1&1&0&0\\0&1&0&0}.
\end{align*}

\begin{figure}[htbp]\centering\begin{tikzpicture}
				\tikzset{enclosed/.style={draw, circle, inner sep=0pt, minimum size=.15cm, fill=black}}
				
				\node[enclosed, label={below, yshift=.0cm: d}] (d) at (0, 0.25) {};
				\node[enclosed, label={below, yshift=-0.07cm: c}] (c) at (2, 0.25) {};
				\draw (c) -- (d) node[midway, right] (edge1) {};
				
				\node[enclosed, label={above, yshift=.0cm: a}] (a) at (0, 2.00) {};
				\node[enclosed, label={above, xshift=.0cm: b}] (b) at (2, 2.00) {};
				\draw (a) -- (b) node[midway, right] (edge1) {};
				
				\node[enclosed, label={below, yshift=.0cm: d}] (d) at (3.0, 0.25) {};
				\node[enclosed, label={below, yshift=-0.07cm: c}] (c) at (5.0, 0.25) {};
				\draw (c) -- (d) node[midway, right] (edge1) {};
				
				\node[enclosed, label={above, yshift=.0cm: a}] (a) at (3.0, 2.00) {};
				\node[enclosed, label={above, xshift=.0cm: b}] (b) at (5.0, 2.00) {};
				\draw (a) -- (b) node[midway, right] (edge1) {};
				\draw (a) -- (c) node[midway, right] (edge2) {};
				\draw (a) -- (d) node[midway, right] (edge3) {};
				\draw (b) -- (c) node[midway, right] (edge4) {};
				\draw (b) -- (d) node[midway, right] (edge5) {};
				
				\node[enclosed, label={below, yshift=.0cm: d}] (d) at (6.0, 0.25) {};
				\node[enclosed, label={below, yshift=-0.07cm: c}] (c) at (8.0, 0.25) {};
				\draw (c) -- (d) node[midway, right] (edge1) {};
				
				\node[enclosed, label={above, yshift=.0cm: a}] (a) at (6.0, 2.00) {};
				\node[enclosed, label={above, xshift=.0cm: b}] (b) at (8.0, 2.00) {};
				\draw (a) -- (b) node[midway, right] (edge1) {};
				\draw (b) -- (d) node[midway, right] (edge2) {};
				
				\node[enclosed, label={below, yshift=.0cm: d}] (d) at (9.0, 0.25) {};
				\node[enclosed, label={below, yshift=-0.07cm: c}] (c) at (11.0, 0.25) {};
				
				\node[enclosed, label={above, yshift=.0cm: a}] (a) at (9.0, 2.00) {};
				\node[enclosed, label={above, xshift=.0cm: b}] (b) at (11.0, 2.00) {};
				\draw (a) -- (b) node[midway, right] (edge1) {};
				\draw (a) -- (c) node[midway, right] (edge2) {};
				\draw (b) -- (c) node[midway, right] (edge3) {};
				\draw (b) -- (d) node[midway, right] (edge4) {};		
			\end{tikzpicture}
\caption{The loop graphs $2K_2$, $K_4$, $P_4$, and a triangle with a pendant vertex, in the order of Figure~\ref{en:fig:forbidden}.}\label{en:fig:graphs}\end{figure}

Each matrix has determinant $1$ over $\F$ and therefore has rank $4$. Each corresponding hieroglyph has genus $2$ and one boundary component.
\begin{lemma}\label{en:lem:forbidden}
If $\rank M(H)>2$, then $H$ contains one of the four hieroglyphs in (B).
\end{lemma}
\begin{proof}
Choose $a,b$ with $m_{ab}=1$, and write, with these indices first,
\[
M=\mat{J&B\\B^{\mathsf T}&C},\qquad J=\mat{0&1\\1&0}.
\]
Simultaneous elementary row and column operations give the congruence
$M\sim J\oplus(C+B^{\mathsf T}JB)$. The second summand is symmetric with zero diagonal. If $\rank M>2$, this summand has a nonzero off-diagonal entry indexed by two loops $c,d$. The principal submatrix of $M$ on $a,b,c,d$ has determinant $1$: its Schur complement with respect to $J$ is $J$.

It remains to classify four-letter cyclic words with nonsingular interlacement matrix. For a $4\times4$ interlacement matrix,
\[
\det M=\operatorname{Pf}(M)^2=(m_{ab}m_{cd}+m_{ac}m_{bd}+m_{ad}m_{bc})^2.
\]
Here $\operatorname{Pf}$ is the Pfaffian of an alternating matrix of order $4$. Over $\F$, subtraction is the same as addition and $x^2=x$, so the determinant equals the Pfaffian. This Pfaffian counts perfect matchings of the loop graph modulo $2$: its three terms correspond to the three partitions of four vertices into pairs. Thus the number of perfect matchings must be odd.

A nonsingular matrix has no zero row, so the loop graph has no isolated vertex. If disconnected, it is $2K_2$. If connected with three edges, it is $P_4$ or $K_{1,3}$. With four edges, the two missing edges of $K_4$ are either disjoint or incident, giving $C_4$ or a triangle with a pendant vertex. With five edges it is $K_4$ minus an edge, and with six it is $K_4$. Their numbers of perfect matchings are, respectively, $1,1,0,2,1,2,3$. Thus only the following four graphs have nonsingular adjacency matrices over $\F$:
\begin{center}\begin{tabular}{ll}
Representative&Loop graph\\\hline
$(ababcdcd)$&$2K_2$\\
$(abcdabcd)$&$K_4$\\
$(abacdcbd)$&$P_4$\\
$(abcadbdc)$&triangle with a pendant vertex
\end{tabular}\end{center}
We show that each determines a unique class of four-letter words up to cyclic shift and relabelling, without using reversal.

Fix two interlacing loops $x,y$. Their ends occur as $x,y,x,y$ and divide the circle into four intervals $I_0,I_1,I_2,I_3$, numbered immediately after these ends. For a third loop with ends in distinct intervals, alternation gives the table
\begin{center}\begin{tabular}{ll}
Intervals containing the ends&Interlacement with $x,y$\\\hline
$I_0,I_1$ or $I_2,I_3$&with $y$ only\\
$I_0,I_3$ or $I_1,I_2$&with $x$ only\\
$I_0,I_2$ or $I_1,I_3$&with both $x$ and $y$.
\end{tabular}\end{center}
If both ends lie in one interval, the loop interlaces with neither $x$ nor $y$.

For $2K_2$, choose $x,y$ from one edge. Each remaining loop has both ends in one interval. Since these two loops interlace, their intervals coincide and their ends alternate there. A cyclic shift and relabelling give $(ababcdcd)$.

For $K_4$, cut the circle at the two ends of $a$. Each remaining loop has one end on each of the two arcs. Pairwise interlacement means that these ends occur in the same order on both arcs when the circle is traversed in its fixed direction. Labelling the three loops $b,c,d$ in this order gives $(abcdabcd)$.

For the path $p-x-y-q$, choose the middle edge $xy$. The loop $p$ interlaces with $x$ only. Moving the starting point from the first $x$ to the second, if necessary, places its ends in $I_0,I_3$. The ends of $q$ lie either in $I_0,I_1$ or in $I_2,I_3$. In the first case, the end of $p$ precedes that of $q$ in their common interval $I_0$; the reverse order would create the unwanted edge $pq$. In the second case, the end of $q$ precedes that of $p$ in $I_3$ for the same reason. The resulting words are
\[
(xpqyqxyp),\qquad(xpyxqyqp).
\]
Interchanging $x,y$ and $p,q$ in the first gives the cyclic shift of the second starting at its sixth letter. Reading the first from its last $p$ and relabelling $p\mapsto a$, $x\mapsto b$, $q\mapsto c$, $y\mapsto d$ gives $(abacdcbd)$.

Finally, suppose the graph is a triangle on $a,b,c$ with $d$ adjacent only to $b$. The triangle has order $a,b,c,a,b,c$, by the same two-arc argument as for $K_4$. Applying the table to $a,c$ shows that both ends of $d$ lie in one interval between their ends. To interlace with $b$, this interval must contain an end of $b$, with the ends of $d$ on opposite sides of it. There are exactly two such intervals, each containing one end of $b$. Thus $d$ is inserted immediately before and after one occurrence of $b$. The two choices are related by a cyclic shift through three ends of the triangle, giving the class $(abcadbdc)$. The chosen four loops therefore form one of the required hieroglyphs.
\end{proof}
\begin{claim}[(A)$\Longleftrightarrow$(E)$\Longleftrightarrow$(B)]
Conditions (A), (B), and (E) are equivalent.
\end{claim}
\begin{proof}
By Lemmas~\ref{en:lem:genus} and~\ref{en:lem:rank}, embeddability in the torus is equivalent to $g(H)\le1$, hence to $\rank M(H)\le2$. Deleting loops amounts to passing to a principal submatrix of $M(H)$; its rank cannot exceed that of $M(H)$. Thus (E) excludes the four nonsingular matrices $M_i$. The converse follows from Lemma~\ref{en:lem:forbidden}.
\end{proof}
For loop graphs the relevant subgraphs are \emph{induced}: vertices and all incident edges are deleted, but individual edges are not. For example, $K_4$ contains $2K_2$ as a subgraph but not as an induced subgraph. Hence neither prohibition can be discarded on that basis.
\begin{claim}[(E)$\Longleftrightarrow$(C)]
The rank condition is equivalent to the stated condition on the loop graph.
\end{claim}
\begin{proof}
If $M=0$, all vertices are isolated. Otherwise choose $a,b$ with $m_{ab}=1$. When $\rank M=2$, the Schur complement in the proof of Lemma~\ref{en:lem:forbidden} vanishes. Consequently, for all $i,j$,
\begin{equation}\label{en:eq:types}m_{ij}=m_{ia}m_{jb}+m_{ib}m_{ja}.
\end{equation}
Assign to vertex $i$ the type $t_i=(m_{ia},m_{ib})\in\F^2$. Type $00$ gives isolated vertices. Within each of $10,01,11$ there are no edges, while between any two distinct nonzero types all edges occur. This gives a complete bipartite or tripartite graph and isolated vertices, as in Figure~\ref{en:fig:types}. Each node of the triangle in that figure denotes an entire part, and a line denotes all edges between the two parts.

Conversely, assign distinct nonzero vectors in $\F^2$ to the two or three parts and the zero vector to isolated vertices. Let $T$ have these vectors as columns. Then $M=T^{\mathsf T}JT$, so $\rank M\le2$. The product of the part sizes may be quadratic in $n$; the linear-time algorithm need not construct the loop graph.
\end{proof}

\begin{figure}[htbp]\centering
\begin{tikzpicture}[every node/.style={draw,circle,minimum size=8mm}]
\node (z) at (0,0) {$00$};\node[draw=none,rectangle] at (0,-.85) {isolated vertices};
\node (a) at (3,0) {$10$};\node (b) at (5,0) {$01$};\node (c) at (4,-1.5) {$11$};\draw (a)--(b)--(c)--(a);
\node[draw=none,rectangle,align=left] at (8,-.6) {each node is a part;\\all edges between parts};
\end{tikzpicture}
\caption{The four vector types in $\F^2$. The zero type is isolated; the three nonzero types give the parts of the loop graph.}\label{en:fig:types}\end{figure}

\FloatBarrier
\subsection*{Reductions and genus}
Figure~\ref{en:fig:reductions} relates the word operations to changes of basis for the intersection form. The lower diagram shows which parallel loop is deleted in $(\beta)$; letters in $U,V$ retain their order.
\begin{lemma}\label{en:lem:reductions}
Both $(\alpha)$ and $(\beta)$ preserve the genus of the disc with bands. Each decreases $n$ and $F(H)$ by one.
\end{lemma}
\begin{proof}
In $(\alpha)$ the deleted row and column are zero: no other loop end lies between the adjacent ends of the deleted loop. Their removal preserves rank.

In $(\beta)$, loops $a,b$ do not interlace and have identical interlacements with every other loop, since no end of a third loop separates either adjacent pair. Thus rows $a,b$ coincide. Replace $e_b$ by $e_a+e_b$ in the homology basis. This vector lies in the radical, and
\[
M(H)\sim M(H')\oplus(0),
\]
where $H'$ is obtained by deleting $b$. The same congruence holds for $(\alpha)$. By Lemma~\ref{en:lem:rank}, the genus is unchanged, and~\eqref{en:eq:euler} then gives $F(H')=F(H)-1$. Inverse operations also preserve genus.
\end{proof}

\begin{figure}[htbp]\centering
\begin{tabular}{lll}
$(\alpha)$&$(aaU)\longrightarrow(U)$&$M\sim M'\oplus(0)$\\[3pt]
$(\beta)$&$(abUbaV)\longrightarrow(aUaV)$&$e_b\longmapsto e_a+e_b$
\end{tabular}\par\medskip
\begin{tikzpicture}
\draw (0,0) ellipse (1.7 and .65);\draw (5,0) ellipse (1.7 and .65);
\foreach \x/\s in {-1.1/a,-.5/b,.5/b,1.1/a}{\node[fill=white] at (\x,.5) {$\s$};}
\node[fill=white] at (0,.65) {$U$};\node[fill=white] at (0,-.65) {$V$};
\node[fill=white] at (4,.5) {$a$};\node[fill=white] at (6,.5) {$a$};\node[fill=white] at (5,.65) {$U$};\node[fill=white] at (5,-.65) {$V$};\draw[->] (2,0)--(3,0);
\end{tikzpicture}
\caption{Word reductions and their effect on the intersection form. The lower diagram depicts $(\beta)$; the rest of the circle runs through $V$.}\label{en:fig:reductions}\end{figure}

\begin{lemma}\label{en:lem:terminal}
If $g(H)\le1$ and neither $(\alpha)$ nor $(\beta)$ applies, then $H$ is $()$, $(abab)$ or $(abcabc)$ up to cyclic shift and relabelling.
\end{lemma}
\begin{proof}
Assume $n>0$ and use the cellular embedding in $\widehat S(H)$. The degree of a face is the number of band sides in its boundary walk, counted with multiplicity. The sum of face degrees is $2n$. We explicitly relate small faces to reductions. Number the $2n$ ends, let $\tau$ exchange the two ends of each loop, and let $\sigma$ advance to the next end in the cyclic order. A boundary step goes along a band side from $i$ to $\tau(i)$ and then along the disc boundary to $\sigma(\tau(i))$. Faces of degree $d$ thus correspond to cycles of length $d$ of $\sigma\tau$.

A degree-one face means $\sigma(\tau(i))=i$, so the two ends $\tau(i),i$ are adjacent and $(\alpha)$ applies. Conversely, adjacent occurrences $aa$ give such a fixed point of $\sigma\tau$ and hence a degree-one face.

For a degree-two face there are distinct ends $i,j$ with
\[
\sigma(\tau(i))=j,\qquad\sigma(\tau(j))=i.
\]
They cannot belong to one loop: $j=\tau(i)$ would imply $\sigma(j)=j$, whereas a cyclic permutation of $2n\ge2$ ends has no fixed point. Thus $i,\tau(i),j,\tau(j)$ are four distinct ends. Call the loops $a,b$. The equations specify adjacent pairs $\tau(i),j$ and $\tau(j),i$, giving $ab$ and $ba$. The word is therefore $(abUbaV)$ and $(\beta)$ applies. Conversely, these adjacent pairs give a cycle of length two and a degree-two face.

Since neither reduction applies, all faces have degree at least three, so
\[
2n\ge3F(H)=3(n+1-2g(H)).
\]
This is impossible for $g(H)=0$; for $g(H)=1$ it gives $n\le3$. If $n=1$, operation $(\alpha)$ applies. If $n=2$, excluding $(\alpha)$ leaves only $(abab)$. If $n=3$, genus one guarantees two interlacing loops $a,b$. Their ends divide the word into four intervals. Both ends of the third loop $c$ in one interval would allow $(\alpha)$. Of the six pairs of distinct intervals, the four adjacent pairs allow $(\beta)$ with $a$ or $b$, while the two opposite pairs give $(abcabc)$ up to relabelling and cyclic shift. These exhaust the possibilities.
\end{proof}
\begin{claim}[(C)$\Longrightarrow$(D)$\Longrightarrow$(A)]
The remaining implications in Theorem~\ref{en:thm:main} hold.
\end{claim}
\begin{proof}
Condition (C) implies (E), so $g(H)\le1$. Apply any available reduction until none remains. Each deletes a loop, so the process terminates. By Lemma~\ref{en:lem:reductions}, the genus remains at most one. By Lemma~\ref{en:lem:terminal}, the terminal hieroglyph is one of those listed in (D). This proves (C)$\Longrightarrow$(D).

Conversely, $()$, $(abab)$, $(abcabc)$ have ranks $0,2,2$ and genera $0,1,1$. Reversing the reductions and using Lemma~\ref{en:lem:reductions} gives $g(H)\le1$, hence embeddability in the torus. Figure~\ref{en:fig:torus} displays the three terminal forms on a square model of the torus, with consistently oriented identifications of opposite sides.

We check that these drawings give the required embeddings after the sides are identified. Identify each square with $[0,1]^2$, let $q:\mathbb R^2\to\mathbb R^2/\mathbb Z^2$ be the quotient map, and put the common vertex at $v=(1/2,1/2)$. In the middle square the loops are $q(v+t(0,1))$ and $q(v+t(1,0))$, $0\le t\le1$. In the right square $a,b,c$ correspond to $(0,1),(1,1),(1,0)$. In particular, the two diagonal halves of $b$ join through identified corners.

Each vector is primitive, so its loop has no self-intersections. The determinant of any pair of distinct vectors is $1$ or $-1$. Hence if $q(v+tu)=q(v+sw)$ for $0\le t,s<1$, then $tu-sw\in\mathbb Z^2$, and invertibility of the integer matrix $(u,-w)$ over $\mathbb Z$ gives $t=s=0$. Distinct loops therefore meet only at their common vertex. Clockwise around the right drawing, the rays occur as $a,b,c,a,b,c$; the middle drawing gives $a,b,a,b$. Choosing the corresponding orientation of the torus realizes the prescribed orders. The left drawing is a single vertex with no loops.
\end{proof}

\begin{figure}[htbp]\centering\begin{tikzpicture}
\foreach \x in {0,3,6}{
  \draw (\x,0) rectangle (\x+2,2);
  \fill (\x+1,1) circle (.125);
}
\draw[very thick] (4,0)--(4,2) (3,1)--(5,1);
\draw[very thick] (7,0)--(7,2) (6,1)--(8,1) (6,0)--(8,2);
\foreach \x in {0,3,6}{
  \foreach \y in {0,2}{
    \foreach \d in {-.10,.10}{
      \fill (\x+1+\d-.05,\y-.055)
        --(\x+1+\d+.05,\y)
        --(\x+1+\d-.05,\y+.055)--cycle;
    }
  }
  \foreach \s in {0,2}{
    \fill (\x+\s-.055,1-.05)
      --(\x+\s+.055,1-.05)
      --(\x+\s,1+.05)--cycle;
  }
}
\node[above] at (4,2) {$a$};
\node[right] at (5,1) {$b$};
\node[above] at (7,2) {$a$};
\node[right] at (8,2) {$b$};
\node[right] at (8,1) {$c$};
\end{tikzpicture}
\caption{Embeddings of $()$, $(abab)$, $(abcabc)$ in the torus represented by a square with opposite sides identified. Loop segments continue through identified sides and corners.}\label{en:fig:torus}\end{figure}

All five conditions of Theorem~\ref{en:thm:main} are now equivalent. The proof also shows that every complete reduction sequence for a hieroglyph of genus at most one ends in one of the three listed forms. No uniqueness of the terminal word is asserted.

\section{Boundary walks and a linear-time algorithm}
For $n>0$, number the occurrences $0,\ldots,2n-1$. Let $\sigma=(0\ 1\ \cdots\ 2n-1)$ be the cyclic-order permutation, and let $\tau$ be the fixed-point-free involution exchanging the two occurrences of each letter. Composition acts from right to left.
\begin{lemma}\label{en:lem:boundary}
The boundary components of the disc with bands are in bijection with the cycles of $\varphi=\sigma\tau$.
\end{lemma}
\begin{proof}
Choose one of the two consistent boundary-traversal conventions. Going along a band side from $i$ to its mate $\tau(i)$, then along the disc boundary to the next end, gives $\sigma(\tau(i))$. Iterating this step traces one boundary component before returning to the initial band side. All $2n$ band sides occur once. The opposite convention gives $\sigma^{-1}\tau$, which has the same number of cycles: it is conjugate to $(\sigma\tau)^{-1}$ by $\tau$.
\end{proof}
Figure~\ref{en:fig:walk} separates the band and disc-boundary portions of a step. For $(abcabc)$, its two cycles are precisely the walks $c_1,c_2$ in Figure~\ref{en:fig:bands}.

\begin{figure}[htbp]\centering
\begin{tabular}{llll}
Word&$\tau$&Cycles of $\sigma\tau$&$(F,g)$\\\hline
$abab$&$(0\ 2)(1\ 3)$&$(0\ 3\ 2\ 1)$&$(1,1)$\\
$abcabc$&$(0\ 3)(1\ 4)(2\ 5)$&$(0\ 4\ 2)(1\ 5\ 3)$&$(2,1)$
\end{tabular}\par\medskip
$\displaystyle i\ \xrightarrow{\text{band}}\ \tau(i)\ \xrightarrow{\text{disc boundary}}\ \sigma(\tau(i))$
\caption{A boundary step and two examples of cycle counting.}\label{en:fig:walk}\end{figure}

Assume the input is an array of length $2n$ whose entries are loop indices $0,\ldots,n-1$. Array access, comparison of indices, and arithmetic on indices take constant time. Complexity is measured in these operations; storage is measured in stored indices and visited flags.
\begin{theorem}
For a hieroglyph with $n$ loops, embeddability in the torus can be recognized in $O(n)$ operations and $O(n)$ space under this input representation and cost model. The algorithm also computes the minimum genus $g(H)$.
\end{theorem}
\begin{proof}
For the empty word, return $F=1$, $g=0$, and an affirmative answer. For $n>0$, first construct the array $\tau$ in one pass: remember the first occurrence of each letter and, at its second occurrence, record both paired positions. Input validity can also be checked in this pass.

Create a visited array of length $2n$ and set $F=0$. For every unvisited position $i$, increase $F$ by one and mark successive positions
\[
j\longmapsto(\tau(j)+1)\bmod 2n
\]
until reaching an already visited position. Since this is a permutation, every position is visited exactly once during these cycle traversals. By Lemma~\ref{en:lem:boundary}, the counter equals $F(H)$. Compute $g=(n+1-F)/2$ and answer affirmatively exactly when $g\le1$. Correctness follows from Lemma~\ref{en:lem:genus}.

Constructing $\tau$ and counting its boundary-permutation cycles require linearly many index and array operations and linear space. Neither the full interlacement matrix nor the explicit loop graph is constructed.
\end{proof}

\clearpage
\endgroup
\begingroup
\selectlanguage{russian}
\renewcommand{\encodingdefault}{TU}
\let\rmdefault\RussianRomanFamily
\normalfont
\renewcommand{\theoremname}{Теорема}
\renewcommand{\lemmaname}{Лемма}
\renewcommand{\claimname}{Утверждение}
\renewcommand{\defnname}{Определение}
\setcounter{section}{0}
\setcounter{subsection}{0}
\setcounter{equation}{0}
\setcounter{figure}{0}
\setcounter{table}{0}
\setcounter{theorem}{0}
\setcounter{lemma}{0}
\setcounter{claim}{0}
\setcounter{defn}{0}
\setcounter{footnote}{0}
\begin{center}
{\small Русская версия}\par\vspace{8pt}
{\LARGE\bfseries Вложимость иероглифов в тор и форма пересечений\par}
\vspace{1.5em}
Тим Березин\footnote{Основная часть работы выполнена на математико-механическом факультете Санкт-Петербургского государственного университета в 2018--2019 гг. Электронная почта: timurovadia@gmail.com.}\par\vspace{1em}
21 октября 2026 г.
\end{center}\vspace{1.5em}
\begin{abstract}
Для иероглифа, то есть одновершинного графа с системой вращения, доказывается эквивалентность вложимости в тор, отсутствия четырёх запрещённых иероглифов, специального вида графа петель, редуцируемости к одному из трёх иероглифов и условия на ранг матрицы перекрещиваний. Доказательство основано на форме пересечений на первых гомологиях диска с ленточками над $\F$ и равенстве её ранга удвоенному роду. Описаны действие редукций и роль краевых окружностей. Подсчёт циклов перестановки краевого обхода даёт линейный алгоритм распознавания при индексированном представлении входа.
\end{abstract}
\section*{Введение}
Пусть $G$ --- граф с одной вершиной и конечным числом $n$ петель. \emph{Иероглифом} назовём граф $G$, снабжённый ориентированным циклическим порядком выходящих из его единственной вершины концов петель. Другое название --- одновершинный граф с системой вращения или одновершинный граф с вращениями.
\begin{defn}
Иероглиф называется \emph{вложимым в тор}, или \emph{реализуемым на торе}, если его петли можно изобразить на торе без самопересечений и взаимных пересечений вне общей вершины так, чтобы при обходе по тору вокруг вершины выходящие из неё отрезки проходились в соответствии с указанным циклическим порядком.
\end{defn}
Таким образом, петли могут иметь общую вершину, но не другие общие точки. Вложение не обязано быть клеточным. Об обобщениях на графы с несколькими вершинами и другие поверхности см.~\cite{ru:MT,ru:EM}. Для произвольного связного графа переход к одновершинному графу осуществляется стягиванием выбранного остовного дерева с переносом системы вращения; одной гомотопической эквивалентности графов для сохранения задачи о вложении недостаточно.

\emph{Алфавитом} будем называть конечное множество, его элементы --- \emph{буквами}. Обозначая петли буквами и обходя вершину в заданном направлении, запишем иероглиф словом длины $2n$, в котором каждая буква встречается дважды. Слова рассматриваются с точностью до переименования букв и циклической перестановки. Направление обхода фиксируется; обращение слова не включается в это отношение эквивалентности. Пустой иероглиф обозначается $()$. Иероглиф \emph{содержит} другой, если последний получается удалением некоторых букв вместе с обоими их вхождениями.

\emph{Графом петель} иероглифа $H$ называется простой граф $H^*$, вершины которого --- петли $H$; две вершины смежны, если соответствующие буквы чередуются, то есть задают $(abab)$. В этом случае петли \emph{перекрещиваются}. Это комбинаторное свойство циклического порядка, а не дополнительная точка пересечения во вложении.

\emph{Редукция иероглифа} --- композиция следующих преобразований:
\begin{align*}
(\alpha)\quad &(aaU)\longmapsto(U),\\
(\beta)\quad &(abUbaV)\longmapsto(aUaV),\qquad a\ne b.
\end{align*}
Здесь $U,V$ --- возможно пустые фрагменты слова, не содержащие буквы $a$ в преобразовании $(\alpha)$ и букв $a,b$ в преобразовании $(\beta)$. Каждая из остальных букв встречается в этих фрагментах суммарно дважды; в случае $(\beta)$ её вхождения могут находиться как в одном фрагменте, так и в разных. В $(\beta)$ две соседние пары используют четыре различных вхождения. Петля из $(\alpha)$ называется \emph{изолированной петлёй}. Изолированная вершина графа петель не обязательно отвечает такой петле: например, в $(abbacc)$ концы $a$ не соседние. Это различие существенно для доказательства редукции.

Граф называется \emph{полным $k$-дольным}, если его вершины разбиты на $k$ независимых долей и любые две вершины разных долей смежны. Пустые доли разрешены. Под объединением ниже понимается дизъюнктное объединение; множество изолированных вершин может быть пустым.
\begin{theorem}\label{ru:thm:main}
Для иероглифа $H$ эквивалентны условия:
\begin{enumerate}
\item[(A)] $H$ вложим в тор;
\item[(B)] $H$ не содержит $(ababcdcd)$, $(abcdabcd)$, $(abacdcbd)$, $(abcadbdc)$;
\item[(C)] граф петель $H^*$ есть объединение изолированных вершин и полного дву- или трёхдольного графа;
\item[(D)] $H$ редуцируется к $()$, $(abab)$ или $(abcabc)$;
\item[(E)] матрица перекрещиваний $M(H)$ над $\F$ имеет ранг не более $2$.
\end{enumerate}
\end{theorem}
Эквивалентность условий (A) и (E) следует из равенства $\rank M(H)=2g(H)$ над $\F$, где $g(H)$ --- род поверхности, полученной заклеиванием краевых окружностей диска с ленточками. Это равенство будет доказано в разделе~\ref{ru:sec:form}. Критерии обсуждаются в~\cite{ru:S}; связь ленточных поверхностей со словами в свободных группах см.~\cite{ru:C}.

\section{Диск с ленточками и форма пересечений}\label{ru:sec:form}
Отметим на ориентированной границе диска $2n$ непересекающихся отрезков в порядке букв слова. Для каждой буквы приклеим к двум её отрезкам ленточку-прямоугольник. Склейки выбираются так, чтобы ориентация диска продолжалась на все ленточки; скрученные ленточки не допускаются. Полученную связную компактную ориентируемую поверхность обозначим $S(H)$ и назовём \emph{диском с $n$ ленточками}. Её компоненты границы называются краевыми окружностями; их число обозначается $F(H)$. Для $()$ полагаем $S(())=D^2$ и $F(())=1$. Далее иероглиф и соответствующий диск с ленточками могут обозначать один и тот же объект. На рис.~\ref{ru:fig:bands} показаны диск с ленточками для $(abcabc)$ и два его краевых обхода. Для чтения правой части рисунка вхождения букв пронумерованы от $0$ до $5$ в порядке $a,b,c,a,b,c$.

\begin{figure}[htbp]\centering
\begin{minipage}[c]{.22\linewidth}\centering
\includegraphics[width=.95\linewidth]{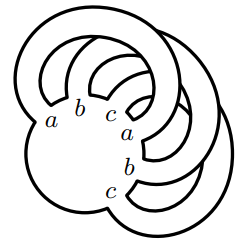}
\end{minipage}\hspace{.06\linewidth}\begin{minipage}[c]{.46\linewidth}\centering
\begin{tikzpicture}[scale=.8,>=Stealth]
\foreach \x/\col/\title/\top/\left/\right in {0/red/{c_1:(0\ 4\ 2)}/{0:a}/{2:c}/{4:b},3.4/blue/{c_2:(1\ 5\ 3)}/{1:b}/{3:a}/{5:c}}{
\begin{scope}[shift={(\x,0)},color=\col]
\draw[->] (90:.75) arc (90:-30:.75);\draw[->] (-30:.75) arc (-30:-150:.75);\draw[->] (-150:.75) arc (-150:-270:.75);
\node at (0,1.55) {$\title$};\node at (0,1.02) {$\top$};\node at (-1.05,-.68) {$\left$};\node at (1.05,-.68) {$\right$};
\end{scope}}
\end{tikzpicture}\end{minipage}
\caption{Иероглиф $(abcabc)$: слева диск с тремя ленточками; справа схема двух краевых окружностей $c_1,c_2$. Метка $i:x$ означает сторону ленточки $x$, отвечающую вхождению $i$; стрелки задают краевой обход.}\label{ru:fig:bands}
\end{figure}

Стягивание диска к вершине, а ленточек к их средним линиям даёт деформационную ретракцию $S(H)$ на граф $H$. В частности,
\[
\chi(S(H))=1-n,\qquad H_1(S(H);\F)\cong\F^n.
\]
Заклеив краевые окружности дисками, получим замкнутую ориентируемую поверхность $\widehat S(H)$; её род обозначим $g(H)$.
\begin{lemma}\label{ru:lem:genus}
Для любого иероглифа
\begin{equation}\label{ru:eq:euler}2g(H)=n+1-F(H).
\end{equation}
Иероглиф вложим в замкнутую ориентируемую поверхность рода $h$ тогда и только тогда, когда $g(H)\le h$.
\end{lemma}
\begin{proof}
Заклеивание $F(H)$ окружностей увеличивает эйлерову характеристику на $F(H)$, откуда $2-2g(H)=1-n+F(H)$. Если $H$ вложен в поверхность рода $h$, его замкнутая регулярная окрестность гомеоморфна $S(H)$. Пусть замыкание дополнения имеет $k$ компонент родов $h_1,\ldots,h_k$ и суммарно $F(H)$ краевых окружностей. Аддитивность $\chi$ при склейке по окружностям даёт
\[
h=g(H)+\sum_{i=1}^k h_i+F(H)-k\ge g(H),
\]
поскольку каждая компонента дополнения имеет непустую границу. Обратно, вложим $H$ в $\widehat S(H)$ и при необходимости добавим $h-g(H)$ ручек внутри диска дополнения. Система вращения сохраняется.
\end{proof}
Таким образом, число компонент дополнения произвольного вложения не отождествляется с $F(H)$: одна компонента может иметь несколько краевых окружностей. Равенство~\eqref{ru:eq:euler} применяется к канонической заклеенной поверхности, где вложение клеточно.
\begin{defn}
\emph{Матрица перекрещиваний} $M(H)=(m_{ij})_{i,j=1}^n$ над $\F$ имеет нулевую диагональ, а $m_{ij}=1$ тогда и только тогда, когда петли $i,j$ перекрещиваются. Это матрица смежности графа петель над $\F$.
\end{defn}
Для формы пересечений $\lambda$ на $H_1(S(H);\F)$ её \emph{радикалом} называется подпространство
\[
\rad\lambda=\{x:\lambda(x,y)=0\text{ для всех }y\in H_1(S(H);\F)\}.
\]
Радикал формы в базисе петель отождествляется с ядром матрицы $M(H)$; далее используется также обозначение $\rad M(H)$.
\begin{lemma}\label{ru:lem:rank}
Матрица $M(H)$ представляет форму пересечений на $H_1(S(H);\F)$ в базисе петель. Её радикал порождён классами краевых окружностей с единственным линейным соотношением --- их суммой, равной нулю. В частности,
\begin{equation}\label{ru:eq:rank}
\dim\rad M(H)=F(H)-1,\qquad \rank M(H)=2g(H).
\end{equation}
\end{lemma}
\begin{proof}
Представим каждый базисный цикл средней линией его ленточки и дугой в диске, соединяющей её концы. Дуги можно выбрать в общем положении. Число их пересечений по модулю $2$ равно $1$ ровно при чередовании концов на границе диска; вне диска циклы разных ленточек не пересекаются. Самопересечение равно нулю, поскольку поверхность ориентируема. Это доказывает утверждение о матрице.

Для полноты опишем форму в стандартном базисе ориентируемой поверхности рода $g=g(H)$ с $F=F(H)$ краевыми окружностями. Выберем циклы ручек
\[
u_1,v_1,\ldots,u_g,v_g
\]
и краевые циклы $c_1,\ldots,c_{F-1}$. Разрезание по циклам ручек и дугам между краевыми окружностями даёт стандартный многоугольник и показывает, что эти $2g+F-1=n$ классов образуют базис. Единственные ненулевые попарные пересечения --- $u_i\cdot v_i=1$. Следовательно, матрица формы имеет вид
\[
\bigoplus_{i=1}^g\mat{0&1\\1&0}\ \oplus\ 0_{F-1}.
\]
Последняя краевая окружность гомологична сумме остальных. Отсюда следуют описание радикала и оба равенства~\eqref{ru:eq:rank}.
\end{proof}
Описание радикала имеет также естественную гомологическую формулировку. Обозначим $S=S(H)$ и $\widehat S=\widehat S(H)$. Для включения $j:S\hookrightarrow\widehat S$ точная последовательность пары и вырезание заклеивающих дисков дают (всюду с коэффициентами в $\F$)
\[
0\longrightarrow H_2(\widehat S)\longrightarrow H_2(\widehat S,S)\longrightarrow H_1(S)\xrightarrow{j_*}H_1(\widehat S)\longrightarrow0.
\]
Здесь первые две ненулевые группы равны $\F$ и $\F^F$; первый ненулевой гомоморфизм переводит $1$ в $(1,\ldots,1)$, а следующий переводит относительные гомологические классы заклеивающих дисков в классы их границ. Следовательно,
\[
H_1(S(H);\F)/\rad M(H)\cong H_1(\widehat S(H);\F).
\]
Это изоморфизм пространств с формой пересечений; справа форма невырождена. Таким образом, радикал отвечает именно краю диска с ленточками, а ненулевая часть формы --- ручкам поверхности. Стандартные сведения о точной последовательности пары и формах на ориентируемых поверхностях см.~\cite{ru:H}.

Рассмотрим вновь $(abcabc)$ на рис.~\ref{ru:fig:bands}. Каждый из двух краевых обходов проходит по одной стороне каждой из трёх ленточек. После ретракции на граф оба обхода поэтому задают один и тот же класс:
\[
[c_1]=[c_2]=e_a+e_b+e_c\quad\text{в }H_1(S(H);\F).
\]
Здесь $e_a,e_b,e_c$ --- базисные классы петель. Действительно,
\[
M(abcabc)=\mat{0&1&1\\1&0&1\\1&1&0},\qquad
M(abcabc)\mat{1\\1\\1}=0.
\]
Радикал одномерен и порождён этим ненулевым классом; единственное независимое соотношение между краевыми классами имеет вид $[c_1]+[c_2]=0$. Способ получения двух циклов, указанных на рисунке, будет обоснован в лемме~\ref{ru:lem:boundary}. Рисунок~\ref{ru:fig:intersection} поясняет вычисление отдельного элемента матрицы формы: чередующиеся концы дают одно пересечение дуг в диске по модулю $2$.

\begin{figure}[htbp]\centering
\begin{tikzpicture}[scale=.75]
\draw (0,0) circle (1);\draw[thick] (-.707,.707)--(.707,-.707);\draw[thick] (.707,.707)--(-.707,-.707);
\node at (-.9,.9) {$a$};\node at (.9,.9) {$b$};\node at (.9,-.9) {$a$};\node at (-.9,-.9) {$b$};
\node[right] at (2,0) {$\displaystyle M(abab)=\mat{0&1\\1&0},\quad\rank M=2,\ g=1,\ F=1.$};
\end{tikzpicture}
\caption{Чередование концов даёт единицу формы пересечений. Пересечение дуг внутри диска служит для вычисления формы на циклах; оно не является пересечением петель во вложении иероглифа.}\label{ru:fig:intersection}
\end{figure}

\section{Доказательство теоремы~\ref{ru:thm:main}}
\subsection*{Запрещённые иероглифы и малые матрицы}
На рис.~\ref{ru:fig:forbidden} приведены четыре запрещённых иероглифа из условия (B). Их графы петель на рис.~\ref{ru:fig:graphs} позволяют прочитать элементы следующих матриц: каждому ребру между двумя различными буквами отвечает единица.

\begin{figure}[htbp]\centering\includegraphics[width=\linewidth]{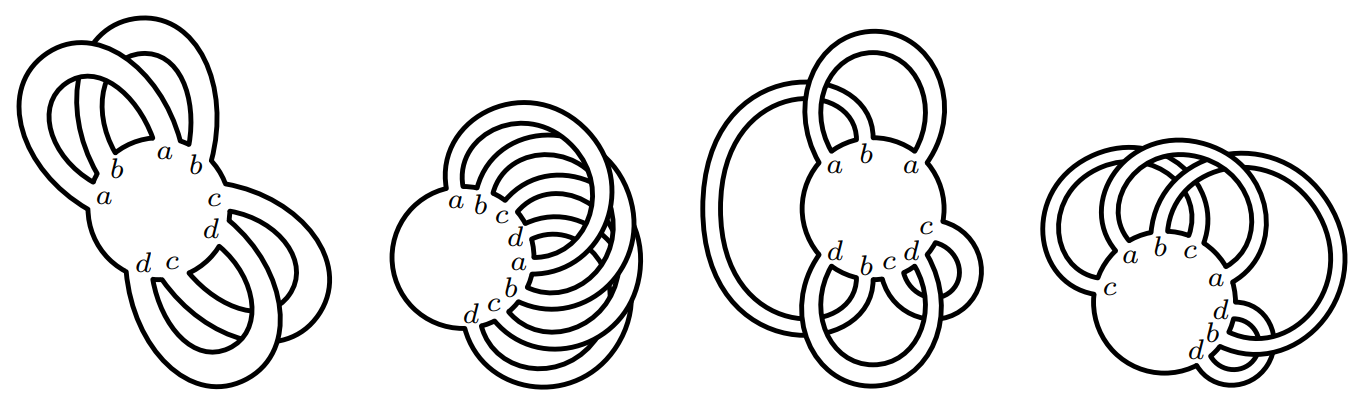}
\caption{Иероглифы $(ababcdcd)$, $(abcdabcd)$, $(abacdcbd)$, $(abcadbdc)$.}\label{ru:fig:forbidden}\end{figure}

Для четырёх иероглифов из (B), в порядке их перечисления и в базисе $a,b,c,d$, матрицы перекрещиваний равны
\begin{align*}
M_1&=\mat{0&1&0&0\\1&0&0&0\\0&0&0&1\\0&0&1&0},&
M_2&=\mat{0&1&1&1\\1&0&1&1\\1&1&0&1\\1&1&1&0},\\
M_3&=\mat{0&1&0&0\\1&0&0&1\\0&0&0&1\\0&1&1&0},&
M_4&=\mat{0&1&1&0\\1&0&1&1\\1&1&0&0\\0&1&0&0}.
\end{align*}

\begin{figure}[htbp]\centering\begin{tikzpicture}
				\tikzset{enclosed/.style={draw, circle, inner sep=0pt, minimum size=.15cm, fill=black}}
				
				\node[enclosed, label={below, yshift=.0cm: d}] (d) at (0, 0.25) {};
				\node[enclosed, label={below, yshift=-0.07cm: c}] (c) at (2, 0.25) {};
				\draw (c) -- (d) node[midway, right] (edge1) {};
				
				\node[enclosed, label={above, yshift=.0cm: a}] (a) at (0, 2.00) {};
				\node[enclosed, label={above, xshift=.0cm: b}] (b) at (2, 2.00) {};
				\draw (a) -- (b) node[midway, right] (edge1) {};
				
				\node[enclosed, label={below, yshift=.0cm: d}] (d) at (3.0, 0.25) {};
				\node[enclosed, label={below, yshift=-0.07cm: c}] (c) at (5.0, 0.25) {};
				\draw (c) -- (d) node[midway, right] (edge1) {};
				
				\node[enclosed, label={above, yshift=.0cm: a}] (a) at (3.0, 2.00) {};
				\node[enclosed, label={above, xshift=.0cm: b}] (b) at (5.0, 2.00) {};
				\draw (a) -- (b) node[midway, right] (edge1) {};
				\draw (a) -- (c) node[midway, right] (edge2) {};
				\draw (a) -- (d) node[midway, right] (edge3) {};
				\draw (b) -- (c) node[midway, right] (edge4) {};
				\draw (b) -- (d) node[midway, right] (edge5) {};
				
				\node[enclosed, label={below, yshift=.0cm: d}] (d) at (6.0, 0.25) {};
				\node[enclosed, label={below, yshift=-0.07cm: c}] (c) at (8.0, 0.25) {};
				\draw (c) -- (d) node[midway, right] (edge1) {};
				
				\node[enclosed, label={above, yshift=.0cm: a}] (a) at (6.0, 2.00) {};
				\node[enclosed, label={above, xshift=.0cm: b}] (b) at (8.0, 2.00) {};
				\draw (a) -- (b) node[midway, right] (edge1) {};
				\draw (b) -- (d) node[midway, right] (edge2) {};
				
				\node[enclosed, label={below, yshift=.0cm: d}] (d) at (9.0, 0.25) {};
				\node[enclosed, label={below, yshift=-0.07cm: c}] (c) at (11.0, 0.25) {};
				
				\node[enclosed, label={above, yshift=.0cm: a}] (a) at (9.0, 2.00) {};
				\node[enclosed, label={above, xshift=.0cm: b}] (b) at (11.0, 2.00) {};
				\draw (a) -- (b) node[midway, right] (edge1) {};
				\draw (a) -- (c) node[midway, right] (edge2) {};
				\draw (b) -- (c) node[midway, right] (edge3) {};
				\draw (b) -- (d) node[midway, right] (edge4) {};		
			\end{tikzpicture}
\caption{Графы петель иероглифов $(ababcdcd)$, $(abcdabcd)$, $(abacdcbd)$, $(abcadbdc)$ (слева направо).}\label{ru:fig:graphs}\end{figure}

Их графы петель --- соответственно $2K_2$, $K_4$, путь $P_4$ и треугольник с висячей вершиной. У каждой из этих матриц определитель над $\F$ равен $1$, а ранг равен $4$. Поэтому каждый соответствующий иероглиф имеет род $2$ и одну краевую окружность.
\begin{lemma}\label{ru:lem:forbidden}
Если $\rank M(H)>2$, то $H$ содержит один из четырёх иероглифов из (B).
\end{lemma}
\begin{proof}
Выберем $a,b$ с $m_{ab}=1$. В этом порядке запишем
\[
M=\mat{J&B\\B^{\mathsf T}&C},\qquad J=\mat{0&1\\1&0}.
\]
Одновременные элементарные преобразования строк и столбцов дают $M\sim J\oplus(C+B^{\mathsf T}JB)$. Матрица $C+B^{\mathsf T}JB$ симметрична и имеет нулевую диагональ. Если ранг $M$ больше $2$, то в этом слагаемом есть ненулевой внедиагональный элемент, отвечающий двум петлям $c,d$. Соответствующая главная подматрица $M$ на $a,b,c,d$ имеет определитель $1$: её дополнение Шура к $J$ равно $J$.

Остаётся классифицировать четырёхбуквенные циклические слова с невырожденной матрицей перекрещиваний. Для такой матрицы
\[
\det M=\operatorname{Pf}(M)^2=(m_{ab}m_{cd}+m_{ac}m_{bd}+m_{ad}m_{bc})^2.
\]
Здесь $\operatorname{Pf}$ --- пфаффиан альтернирующей матрицы порядка $4$. Над $\F$ знаки минус совпадают со знаками плюс, а $x^2=x$ для каждого $x\in\F$, поэтому определитель равен самому пфаффиану. Пфаффиан здесь равен числу совершенных паросочетаний графа петель по модулю $2$: каждое из трёх слагаемых отвечает одному разбиению четырёх вершин на две пары. Поэтому требуется нечётность этого числа.

В невырожденной матрице нет нулевой строки, так что граф петель не имеет изолированных вершин. Если он несвязен, это $2K_2$. Если он связен и имеет три ребра, то это путь $P_4$ или звезда $K_{1,3}$. Если рёбер четыре, то два отсутствующих ребра графа $K_4$ либо не имеют общей вершины, либо имеют её; получаются соответственно цикл $C_4$ и треугольник с висячей вершиной. При пяти рёбрах получается $K_4$ без одного ребра, при шести --- $K_4$. Число совершенных паросочетаний в этих графах равно соответственно $1$ для $2K_2$, $1$ для $P_4$, $0$ для $K_{1,3}$, $2$ для $C_4$, $1$ для треугольника с висячей вершиной, $2$ для $K_4$ без ребра и $3$ для $K_4$. Следовательно, невырожденные матрицы смежности над $\F$ имеют только четыре графа из следующей таблицы:
\begin{center}\begin{tabular}{ll}
Представитель&Граф петель\\\hline
$(ababcdcd)$&$2K_2$\\
$(abcdabcd)$&$K_4$\\
$(abacdcbd)$&$P_4$\\
$(abcadbdc)$&треугольник с висячей вершиной
\end{tabular}\end{center}
Докажем, что каждый из этих графов определяет единственный класс четырёхбуквенных слов с точностью до циклического сдвига и переименования букв. Обращение слова не используется.

Выберем две перекрещивающиеся петли $x,y$. Их концы идут в порядке $x,y,x,y$ и делят окружность на четыре промежутка $I_0,I_1,I_2,I_3$, нумеруемые сразу после указанных концов. Для третьей петли, концы которой лежат в разных промежутках, непосредственно из чередования концов получаем:
\begin{center}\begin{tabular}{ll}
Промежутки концов&Перекрещивания с $x,y$\\\hline
$I_0,I_1$ или $I_2,I_3$&только с $y$\\
$I_0,I_3$ или $I_1,I_2$&только с $x$\\
$I_0,I_2$ или $I_1,I_3$&с $x$ и с $y$.
\end{tabular}\end{center}
При обоих концах в одном промежутке перекрещиваний с $x,y$ нет.

Если граф равен $2K_2$, выберем $x,y$ из одного ребра. Каждая из двух оставшихся петель имеет оба конца в одном промежутке. Поскольку эти две петли перекрещиваются, промежуток у них общий. В нём их концы чередуются. Циклический сдвиг и переименование дают $(ababcdcd)$.

Если граф равен $K_4$, разрежем окружность в двух концах петли $a$. Каждая из остальных петель имеет по одному концу на каждой из двух полученных дуг. Попарное перекрещивание остальных петель означает, что их концы идут на обеих дугах в одинаковом порядке при обходе окружности в фиксированном направлении. Обозначая этот порядок через $b,c,d$, получаем $(abcdabcd)$.

Если граф есть путь $p-x-y-q$, выберем его среднее ребро $xy$. Петля $p$ перекрещивается только с $x$ из пары $x,y$. Сдвиг начала обхода от первого $x$ ко второму позволяет считать, что концы $p$ лежат в $I_0,I_3$. По таблице концы $q$ лежат либо в $I_0,I_1$, либо в $I_2,I_3$. В первом случае в общем промежутке $I_0$ конец $p$ предшествует концу $q$: обратный порядок дал бы запрещённое ребро $pq$. Во втором случае в общем промежутке $I_3$ конец $q$ предшествует концу $p$ по той же причине. Получаются соответственно слова
\[
(xpqyqxyp),\qquad(xpyxqyqp).
\]
Замена $x\leftrightarrow y$, $p\leftrightarrow q$ в первом слове даёт циклический сдвиг второго, начинающийся с его шестой буквы. Первое слово, прочитанное с последнего $p$, при переименовании $p\mapsto a$, $x\mapsto b$, $q\mapsto c$, $y\mapsto d$ становится $(abacdcbd)$.

Наконец, пусть граф состоит из треугольника на $a,b,c$ и вершины $d$, смежной только с $b$. Треугольнику соответствует порядок $a,b,c,a,b,c$: это тот же аргумент с двумя дугами, что и для $K_4$. Применяя таблицу к перекрещивающейся паре $a,c$, видим, что оба конца $d$ лежат в одном из четырёх промежутков между концами $a,c$. Чтобы $d$ перекрещивалась с $b$, этот промежуток должен содержать один конец $b$, а концы $d$ должны лежать по разные стороны от него. В порядке $a,b,c,a,b,c$ имеются ровно два таких промежутка, каждый с единственным концом $b$. Таким образом, $d$ вставляется непосредственно перед одним вхождением $b$ и непосредственно после него. Два выбора вхождения $b$ переходят друг в друга циклическим сдвигом на три конца треугольника. Получается класс $(abcadbdc)$. Значит, выбранные четыре петли задают один из требуемых иероглифов.
\end{proof}
\begin{claim}[(A)$\Longleftrightarrow$(E)$\Longleftrightarrow$(B)]
Условия (A), (B) и (E) эквивалентны.
\end{claim}
\begin{proof}
По леммам~\ref{ru:lem:genus},~\ref{ru:lem:rank} вложимость в тор равносильна $g(H)\le1$, то есть $\rank M(H)\le2$. Удаление петель означает переход к главной подматрице $M(H)$, ранг которой не превосходит ранга $M(H)$. Поэтому (E) исключает четыре невырожденные матрицы $M_i$. Обратная импликация следует из леммы~\ref{ru:lem:forbidden}.
\end{proof}
При переходе к графам петель здесь необходимо говорить об \emph{индуцированных} подграфах: удаляются вершины вместе со всеми инцидентными рёбрами, но не отдельные рёбра. В частности, $K_4$ содержит $2K_2$ как обычный подграф, однако не как индуцированный; поэтому эти запреты нельзя отбрасывать один на основании другого.
\begin{claim}[(E)$\Longleftrightarrow$(C)]
Условие на ранг равносильно условию на граф петель.
\end{claim}
\begin{proof}
Если $M=0$, граф состоит из изолированных вершин. Иначе выберем $a,b$ с $m_{ab}=1$. При $\rank M=2$ дополнение Шура из доказательства леммы~\ref{ru:lem:forbidden} равно нулю. Поэтому для любых $i,j$
\begin{equation}\label{ru:eq:types}m_{ij}=m_{ia}m_{jb}+m_{ib}m_{ja}.
\end{equation}
Сопоставим вершине $i$ тип $t_i=(m_{ia},m_{ib})\in\F^2$. Тип $00$ даёт изолированные вершины. Внутри каждого из типов $10,01,11$ рёбер нет, а между двумя различными ненулевыми типами все рёбра есть. Получается полный дву- или трёхдольный граф и изолированные вершины, как показано на рис.~\ref{ru:fig:types}. Каждый из трёх соединённых узлов рисунка обозначает целую долю; соединяющая линия означает наличие всех рёбер между двумя долями.

Обратно, присвоим двум или трём долям различные ненулевые векторы из $\F^2$, а изолированным вершинам вектор $0$. Пусть $T$ --- матрица с этими столбцами. Тогда $M=T^{\mathsf T}JT$, так что $\rank M\le2$. В частности, произведение размеров долей может быть квадратичным по $n$; строить весь граф петель для линейного алгоритма не требуется.
\end{proof}

\begin{figure}[htbp]\centering
\begin{tikzpicture}[every node/.style={draw,circle,minimum size=8mm}]
\node (z) at (0,0) {$00$};\node[draw=none,rectangle] at (0,-.85) {\shortstack{изолированные\\вершины}};
\node (a) at (3,0) {$10$};\node (b) at (5,0) {$01$};\node (c) at (4,-1.5) {$11$};\draw (a)--(b)--(c)--(a);
\node[draw=none,rectangle,align=left] at (8.2,-.6) {каждый узел --- доля;\\между долями все рёбра};
\end{tikzpicture}
\caption{Четыре типа векторов в $\F^2$: нулевой тип изолирован, три ненулевых типа образуют доли графа петель.}\label{ru:fig:types}\end{figure}

\FloatBarrier
\subsection*{Редукции и род}
На рис.~\ref{ru:fig:reductions} сопоставлены преобразования слова и соответствующие преобразования базиса формы пересечений. Нижняя схема показывает, какая из двух параллельных петель удаляется при $(\beta)$; буквы в фрагментах $U,V$ остаются на месте.
\begin{lemma}\label{ru:lem:reductions}
Оба преобразования $(\alpha)$ и $(\beta)$ сохраняют род диска с ленточками. При каждом из них $n$ и $F(H)$ уменьшаются на единицу.
\end{lemma}
\begin{proof}
В случае $(\alpha)$ строка и столбец удаляемой петли нулевые: между её соседними концами нет конца другой петли. Удаление этой строки и столбца сохраняет ранг.

В случае $(\beta)$ петли $a,b$ не перекрещиваются и имеют одинаковые перекрещивания со всеми остальными петлями. Действительно, две пары соседних концов не разделяются концами третьей петли. Значит, строки $a,b$ матрицы совпадают. В базисе гомологий заменим $e_b$ на $e_a+e_b$: этот вектор лежит в радикале формы, и
\[
M(H)\sim M(H')\oplus(0),
\]
где $H'$ получается удалением $b$. Такое же равенство верно для $(\alpha)$. По лемме~\ref{ru:lem:rank} род сохраняется; затем~\eqref{ru:eq:euler} даёт $F(H')=F(H)-1$. Обратные операции также сохраняют род.
\end{proof}

\begin{figure}[htbp]\centering
\begin{tabular}{lll}
$(\alpha)$&$(aaU)\longrightarrow(U)$&$M\sim M'\oplus(0)$\\[3pt]
$(\beta)$&$(abUbaV)\longrightarrow(aUaV)$&$e_b\longmapsto e_a+e_b$
\end{tabular}\par\medskip
\begin{tikzpicture}
\draw (0,0) ellipse (1.7 and .65);\draw (5,0) ellipse (1.7 and .65);
\foreach \x/\s in {-1.1/a,-.5/b,.5/b,1.1/a}{\node[fill=white] at (\x,.5) {$\s$};}
\node[fill=white] at (0,.65) {$U$};\node[fill=white] at (0,-.65) {$V$};
\node[fill=white] at (4,.5) {$a$};\node[fill=white] at (6,.5) {$a$};\node[fill=white] at (5,.65) {$U$};\node[fill=white] at (5,-.65) {$V$};\draw[->] (2,0)--(3,0);
\end{tikzpicture}
\caption{Редукции в циклическом слове и действие на форме пересечений. Нижняя схема показывает $(\beta)$; оставшаяся часть окружности проходит через $V$.}\label{ru:fig:reductions}\end{figure}

\begin{lemma}\label{ru:lem:terminal}
Если $g(H)\le1$ и к $H$ неприменимы $(\alpha)$, $(\beta)$, то $H$ есть $()$, $(abab)$ или $(abcabc)$ с точностью до циклической перестановки и переименования букв.
\end{lemma}
\begin{proof}
Пусть $n>0$. Рассмотрим клеточное вложение в $\widehat S(H)$. Степень грани --- число сторон ленточек в соответствующем краевом обходе, с учётом кратности. Сумма степеней граней равна $2n$. Докажем явно соответствие граней малой степени редукциям. Пронумеруем $2n$ концов петель; через $\tau$ обозначим перестановку, меняющую местами концы одной петли, а через $\sigma$ --- переход к следующему концу в циклическом порядке. Один шаг вдоль края сначала проходит по стороне ленточки от $i$ к $\tau(i)$, затем вдоль границы диска к $\sigma(\tau(i))$. Поэтому граням степени $d$ соответствуют циклы длины $d$ перестановки $\sigma\tau$.

Грань степени $1$ означает $\sigma(\tau(i))=i$. Следовательно, концы $\tau(i),i$ одной петли соседние, и применима $(\alpha)$. Обратно, соседние вхождения $aa$ дают такую неподвижную точку перестановки $\sigma\tau$ и тем самым грань степени $1$.

Для грани степени $2$ имеются различные концы $i,j$, для которых
\[
\sigma(\tau(i))=j,\qquad\sigma(\tau(j))=i.
\]
Они не могут принадлежать одной петле: при $j=\tau(i)$ первое равенство дало бы $\sigma(j)=j$, тогда как циклическая перестановка $2n\ge2$ концов не имеет неподвижных точек. Значит, четыре конца $i,\tau(i),j,\tau(j)$ различны. Пусть соответствующие петли называются $a,b$. Равенства задают соседние пары $\tau(i),j$ и $\tau(j),i$, то есть пары $ab$ и $ba$. Поэтому циклическое слово имеет вид $(abUbaV)$, и применима $(\beta)$. Обратно, две пары соседних концов в слове $(abUbaV)$ дают указанный цикл длины $2$ и грань степени $2$.

Поскольку обе редукции неприменимы, все грани имеют степень не менее $3$:
\[
2n\ge3F(H)=3(n+1-2g(H)).
\]
При $g(H)=0$ это невозможно, а при $g(H)=1$ получаем $n\le3$. При $n=1$ применима $(\alpha)$. При $n=2$ отсутствие $(\alpha)$ оставляет только $(abab)$. При $n=3$ род $1$ означает, что есть две перекрещивающиеся петли $a,b$. Их концы делят циклическое слово на четыре промежутка. Оба конца третьей петли $c$ в одном промежутке дали бы $(\alpha)$. Из шести пар различных промежутков четыре пары соседних промежутков дают $(\beta)$ с $a$ или $b$, а две пары противоположных промежутков дают $(abcabc)$ с точностью до переименования и циклической перестановки. Это исчерпывает случаи.
\end{proof}
\begin{claim}[(C)$\Longrightarrow$(D)$\Longrightarrow$(A)]
Оставшиеся импликации теоремы~\ref{ru:thm:main} справедливы.
\end{claim}
\begin{proof}
Из (C) уже следует (E), а значит, $g(H)\le1$. Выполняем любую доступную редукцию, пока это возможно. Каждая удаляет одну петлю, поэтому процесс конечен. По лемме~\ref{ru:lem:reductions} род остаётся не более $1$. Лемма~\ref{ru:lem:terminal} показывает, что окончательный иероглиф входит в список из (D). Таким образом, требуемая редукция существует для любого иероглифа, удовлетворяющего условию (C).

Обратно, у $()$, $(abab)$, $(abcabc)$ ранги матриц равны соответственно $0,2,2$, а роды --- $0,1,1$. Обращая последовательность редукций, по лемме~\ref{ru:lem:reductions} получаем $g(H)\le1$ и вложимость в тор. На рис.~\ref{ru:fig:torus} показаны вложения трёх конечных форм на развёртке тора; направления отождествляемых сторон квадратов согласованы попарно.

Проверим геометрически эти изображения с учётом склеек. Отождествим каждый квадрат с $[0,1]^2$ и обозначим через $q:\mathbb R^2\to\mathbb R^2/\mathbb Z^2$ проекцию на тор, а через $v=(1/2,1/2)$ общую вершину. На среднем квадрате петли задаются кривыми $q(v+t(0,1))$ и $q(v+t(1,0))$ при $0\le t\le1$. На правом квадрате петлям $a,b,c$ соответствуют векторы $(0,1),(1,1),(1,0)$. В частности, две диагональные половины петли $b$ соединяются через отождествлённые углы квадрата.

Каждый из этих целочисленных векторов примитивен, поэтому соответствующая петля не имеет самопересечений. Определитель любой пары различных векторов равен $1$ или $-1$. Следовательно, если $q(v+tu)=q(v+sw)$ при $0\le t,s<1$, то $tu-sw\in\mathbb Z^2$, а обратимость целочисленной матрицы $(u,-w)$ над $\mathbb Z$ даёт $t=s=0$. Таким образом, разные петли встречаются только в общей вершине. При обходе правого изображения по часовой стрелке лучи идут в порядке $a,b,c,a,b,c$; среднее изображение даёт $a,b,a,b$. Выбирая соответствующую ориентацию тора, получаем заданные циклические порядки. Левое изображение задаёт одну вершину без петель.
\end{proof}

\begin{figure}[htbp]\centering\begin{tikzpicture}
\foreach \x in {0,3,6}{
  \draw (\x,0) rectangle (\x+2,2);
  \fill (\x+1,1) circle (.125);
}
\draw[very thick] (4,0)--(4,2) (3,1)--(5,1);
\draw[very thick] (7,0)--(7,2) (6,1)--(8,1) (6,0)--(8,2);
\foreach \x in {0,3,6}{
  \foreach \y in {0,2}{
    \foreach \d in {-.10,.10}{
      \fill (\x+1+\d-.05,\y-.055)
        --(\x+1+\d+.05,\y)
        --(\x+1+\d-.05,\y+.055)--cycle;
    }
  }
  \foreach \s in {0,2}{
    \fill (\x+\s-.055,1-.05)
      --(\x+\s+.055,1-.05)
      --(\x+\s,1+.05)--cycle;
  }
}
\node[above] at (4,2) {$a$};
\node[right] at (5,1) {$b$};
\node[above] at (7,2) {$a$};
\node[right] at (8,2) {$b$};
\node[right] at (8,1) {$c$};
\end{tikzpicture}
\caption{Вложения иероглифов $()$, $(abab)$, $(abcabc)$ в тор, представленный квадратом с отождествлёнными противоположными сторонами. Отрезки петель продолжаются через склеиваемые стороны и углы.}\label{ru:fig:torus}\end{figure}

Таким образом, все пять условий теоремы~\ref{ru:thm:main} эквивалентны. Доказательство показывает также, что любая полная последовательность редукций иероглифа рода не более $1$ заканчивается одним из трёх указанных иероглифов. Единственность конечного слова здесь не утверждается.

\section{Краевые обходы и линейный алгоритм}
Для $n>0$ пронумеруем вхождения букв числами $0,\ldots,2n-1$. Перестановку циклического порядка обозначим
\[
\sigma=(0\ 1\ \cdots\ 2n-1),
\]
а инволюцию без неподвижных точек, меняющую местами два вхождения каждой буквы, --- через $\tau$. Композиции действуют справа налево.
\begin{lemma}\label{ru:lem:boundary}
Краевые окружности диска с ленточками находятся во взаимно однозначном соответствии с циклами перестановки $\varphi=\sigma\tau$.
\end{lemma}
\begin{proof}
Выберем одну из двух согласованных конвенций обхода ориентированной границы. Перейдя по стороне ленточки от конца $i$ к парному концу $\tau(i)$, затем вдоль границы диска к следующему концу, получим $\sigma(\tau(i))$. Повторение этой операции замыкается ровно после одного краевого обхода. Все $2n$ сторон ленточек участвуют ровно по одному разу. При противоположной конвенции получается $\sigma^{-1}\tau$; число циклов то же, поскольку эта перестановка сопряжена с $(\sigma\tau)^{-1}$ посредством $\tau$.
\end{proof}
На рис.~\ref{ru:fig:walk} один шаг $i\mapsto\sigma(\tau(i))$ разбит на переход по ленточке и переход вдоль границы диска. В строке для $(abcabc)$ получаются ровно два обхода $c_1,c_2$, выделенные на рис.~\ref{ru:fig:bands}.

\begin{figure}[htbp]\centering
\begin{tabular}{llll}
Слово&$\tau$&Циклы $\sigma\tau$&$(F,g)$\\\hline
$abab$&$(0\ 2)(1\ 3)$&$(0\ 3\ 2\ 1)$&$(1,1)$\\
$abcabc$&$(0\ 3)(1\ 4)(2\ 5)$&$(0\ 4\ 2)(1\ 5\ 3)$&$(2,1)$
\end{tabular}\par\medskip
$\displaystyle i\ \xrightarrow{\text{ленточка}}\ \tau(i)\ \xrightarrow{\text{граница диска}}\ \sigma(\tau(i))$
\caption{Один шаг краевого обхода и два примера подсчёта циклов.}\label{ru:fig:walk}\end{figure}

Предполагается, что иероглиф задан массивом длины $2n$, элементы которого --- номера петель от $0$ до $n-1$. Доступ к элементу массива, сравнение индексов и арифметические операции над индексами считаются операциями постоянного времени. Сложность алгоритма оценивается числом таких операций, а используемая память --- числом хранимых индексов и отметок посещения.
\begin{theorem}
Для иероглифа с $n$ петлями существует алгоритм распознавания вложимости в тор, выполняющий $O(n)$ операций и использующий $O(n)$ памяти при указанном представлении входа и способе оценки сложности. Алгоритм одновременно вычисляет минимальный род $g(H)$.
\end{theorem}
\begin{proof}
Для пустого слова возвращаем $F=1$, $g=0$ и положительный ответ. При $n>0$ сначала одним проходом строим массив $\tau$: для каждой буквы запоминаем первое вхождение, а при втором записываем обе парные позиции. При необходимости в том же проходе проверяется корректность входа.

Затем заводим массив отметок посещения длины $2n$ и счётчик $F=0$. Для каждой ещё не отмеченной позиции $i$ увеличиваем $F$ на $1$ и последовательно отмечаем позиции
\[
j\longmapsto(\tau(j)+1)\bmod 2n
\]
до возвращения к уже отмеченной позиции. Поскольку это перестановка, каждая позиция посещается ровно один раз во внутренних обходах. По лемме~\ref{ru:lem:boundary} полученный счётчик равен $F(H)$. Вычисляем $g=(n+1-F)/2$ и отвечаем утвердительно ровно при $g\le1$. Корректность следует из леммы~\ref{ru:lem:genus}.

Построение $\tau$ и подсчёт циклов требуют линейного числа операций над индексами и массивами; память также линейна. Полная матрица перекрещиваний и явный граф петель не строятся.
\end{proof}

\clearpage
\endgroup

\begin{thebibliography}{6}
\small\setlength{\itemsep}{3pt}
\bibitem{en:MT} B. Mohar and C. Thomassen. \emph{Graphs on Surfaces}. Johns Hopkins Studies in the Mathematical Sciences. Johns Hopkins University Press, 2001.
\bibitem{en:EM} J. Ellis-Monaghan and I. Moffatt. Embedded Graphs. In: \emph{Graphs on Surfaces}. SpringerBriefs in Mathematics. Springer, New York, 2013, pp.~1--22. \url{https://doi.org/10.1007/978-1-4614-6971-1_1}.
\bibitem{en:S} A. B. Skopenkov. \emph{Algebraic Topology from a Geometric Viewpoint} (in Russian). MCCME, Moscow, 2020, 302 pp. Updated version: \url{https://www.mccme.ru/circles/oim/obstruct.pdf}.
\bibitem{en:C} M. Culler. Using surfaces to solve equations in free groups. \emph{Topology} \textbf{20} (1981), 133--145. \url{https://core.ac.uk/download/pdf/81186523.pdf} (accessed June 11, 2021).
\bibitem{en:T} C. Thomassen. Embeddings of graphs. \emph{Discrete Mathematics} \textbf{124}(1--3) (1994), 217--228. \url{https://doi.org/10.1016/0012-365X(92)00062-V}.
\bibitem{en:H} A. Hatcher. \emph{Algebraic Topology}. Cambridge University Press, 2002. \url{https://pi.math.cornell.edu/~hatcher/AT/ATpage.html}.
\end{thebibliography}

\begin{thebibliography}{6}
\small\raggedright\setlength{\itemsep}{3pt}
\bibitem{ru:MT} B. Mohar and C. Thomassen. \emph{Graphs on Surfaces}. Johns Hopkins Studies in the Mathematical Sciences. Johns Hopkins University Press, 2001.
\bibitem{ru:EM} J. Ellis-Monaghan and I. Moffatt. Embedded Graphs. In: \emph{Graphs on Surfaces}. SpringerBriefs in Mathematics. Springer, New York, 2013, pp.~1--22. \url{https://doi.org/10.1007/978-1-4614-6971-1_1}.
\bibitem{ru:S} А. Б. Скопенков. \emph{Алгебраическая топология с геометрической точки зрения}. М.: МЦНМО, 2020. 302 с. Обновляемая версия: \url{https://www.mccme.ru/circles/oim/obstruct.pdf}.
\bibitem{ru:C} M. Culler. Using surfaces to solve equations in free groups. \emph{Topology} \textbf{20} (1981), 133--145. \url{https://core.ac.uk/download/pdf/81186523.pdf} (accessed June 11, 2021).
\bibitem{ru:T} C. Thomassen. Embeddings of graphs. \emph{Discrete Mathematics} \textbf{124}(1--3) (1994), 217--228. \url{https://doi.org/10.1016/0012-365X(92)00062-V}.
\bibitem{ru:H} A. Hatcher. \emph{Algebraic Topology}. Cambridge University Press, 2002. \url{https://pi.math.cornell.edu/~hatcher/AT/ATpage.html}.
\end{thebibliography}
\end{document}